\documentclass[11pt]{article}
\usepackage{amssymb,amsthm,amsmath}
\def\qed{\hfill $\Box$}

\newcommand\pf{\smallbreak\noindent \texttt{Proof}. }

\begin{document}

\newtheorem{thm}{Theorem}[section]
\newtheorem{prop}[thm]{Proposition}
\newtheorem{lem}[thm]{Lemma}
\newtheorem{cor}[thm]{Corollary}
\newtheorem{ex}[thm]{Example}
\renewcommand{\thefootnote}{*}

\title{\bf On the automorphism groups of some nilpotent 3-dimensional Leibniz algebras}

\author{\textbf{L.A.~Kurdachenko, O.O.~Pypka}\\
Oles Honchar Dnipro National University, Dnipro, Ukraine\\
{\small e-mail: lkurdachenko@gmail.com, sasha.pypka@gmail.com}\\
\textbf{M.M.~Semko}\\
State Tax University, Irpin, Ukraine\\
{\small e-mail: dr.mykola.semko@gmail.com}}
\date{}

\maketitle

\begin{abstract}
Let $L$ be an algebra over a field $F$ with the binary operations $+$ and $[,]$. Then $L$ is called a left Leibniz algebra if it satisfies the left Leibniz identity: $[[a,b],c]=[a,[b,c]]-[b,[a,c]]$ for all elements $a,b,c\in L$. The structure of the automorphism group of $3$-dimensional Leibniz algebras, which have nilpotency class $2$ and a one-dimensional center, is studied.
\end{abstract}



\thispagestyle{empty}

\section{Introduction.}
Let $L$ be an algebra over a field $F$ with the binary operations $+$ and $[,]$. Then $L$ is called a \textit{left Leibniz algebra} if it satisfies the left Leibniz identity:
$$[[a,b],c]=[a,[b,c]]-[b,[a,c]]$$
for all elements $a,b,c\in L$.

Leibniz algebras appeared first in the paper of A.~Blokh~\cite{BA1965}, but the term ``Leibniz algebra'' appears in the book of J.-L.~Loday~\cite{LJ1992}, and the article of J.-L.~Loday~\cite{LJ1993}. In \cite{LP1993}, J.-L.~Loday and T.~Pirashvili began to study the properties of Leibniz algebras. The theory of Leibniz algebras is developing very intensively in various directions of research, although in some places rather sporadically. Some modern results of this theory can be found in the monograph~\cite{AOR2020}. It is worth noting that Leibniz algebras are a fairly broad generalization of Lie algebras. On the other hand, if $L$ is a Leibniz algebra, in which $[a,a]=0$ for every element $a\in L$, then it is a Lie algebra. Thus, Lie algebras can be characterized as anticommutative Leibniz algebras. At the same time, there is a very significant difference between Lie algebras and Leibniz algebras (see, for example, survey papers \cite{CPSY2019,KKPS2017,KSeSu2020,SI2021}).

Let $L$ be a Leibniz algebra. A linear transformation $f$ of $L$ is called an \textit{endomorphism} of $L$, if
$$f([a,b])=[f(a),f(b)]$$
for all elements $a,b\in L$. Clearly, a product of two endomorphisms of $L$ is also endomorphism, so that the set of all endomorphisms of $L$ is a semigroup by a multiplication. We note that the sum of two endomorphisms of $L$ is not necessary to be an endomorphism of $L$, so that we cannot speak about the endomorphism ring of $L$.

A bijective endomorphism of $L$ is called an \textit{automorphism} of $L$. It is not hard to show that the set $\textrm{Aut}_{[,]}(L)$ of all automorphisms of $L$ is a group by a multiplication (see, for example, \cite{KPS2023}).

As for other algebraic structures, the search for the structure of automorphism groups of Leibniz algebras is one of the important problems of this theory. It should be noted that automorphisms groups of Leibniz algebras have hardly been studied. Among the papers that were devoted to this problem, we can note \cite{AKOZ2020,LRT2016}.

It is natural to start studying automorphism groups of Leibniz algebras, the structure of which has been studied quite fully. A description of the structure of automorphism groups of infinite-dimensional cyclic Leibniz algebras was obtained in \cite{KSuY2022}, and of finite-dimensional cyclic Leibniz algebras was obtained in \cite{KPS2023}. The question naturally arises about automorphism groups of Leibniz algebras of low dimension. The case of two-dimensional Leibniz algebras is quite simple, and the automorphism groups of such algebras were described in \cite{KPV2023}. Unlike Lie algebras, the situation with Leibniz algebras of dimension 3 is very diverse. The most detailed description of three-dimensional Leibniz algebras can be found in \cite{KPS2022}. Due to the large number of types of three-dimensional Leibniz algebras, the study of the structure of their automorphism groups is very voluminous and requires a gradual approach. The first step was taken in the article \cite{KPV2023}, in which the descriptions of the automorphism groups of nilpotent Leibniz algebras with nilpotency class 3 and also of nilpotent Leibniz algebras with nilpotency class 2 and a 2-dimensional center were obtained. In the paper \cite{KPSe2023} the study of the structure of the automorphism groups of nilpotent Leibniz algebras with nilpotency class 2 and 1-dimensional center was started. This article is devoted to the continuation of the study of the automorphism groups of such Leibniz algebras.

\section{Preliminary information about Leibniz algebras.}
Let $L$ be a Leibniz algebra over a field $F$. Then $L$ is called \textit{abelian} if $[a,b]=0$ for every elements $a,b\in L$. In particular, an abelian Leibniz algebra is a Lie algebra.

If $A,B$ are subspaces of $L$, then $[A,B]$ will denote a subspace generated by all elements $[a,b]$ where $a\in A$, $b\in B$. A subspace $A$ of $L$ is called a \textit{subalgebra} of $L$, if $[a,b]\in A$ for every $a,b\in A$. A subalgebra $A$ of $L$ is called a \textit{left} (respectively \textit{right}) \textit{ideal} of $L$, if $[b,a]\in A$ (respectively $[a,b]\in A$) for every $a\in A$, $b\in L$. A subalgebra $A$ of $L$ is called an \textit{ideal} of $L$ (more precisely, \textit{two-sided ideal}) if it is both a left ideal and a right ideal.

Denote by $\mathrm{Leib}(L)$ the subspace generated by the elements $[a,a]$, $a\in L$. It is not hard to prove that $\mathrm{Leib}(L)$ is an ideal of $L$. The ideal $\mathrm{Leib}(L)$ is called the \textit{Leibniz kernel} of $L$.

The \textit{left} (respectively \textit{right}) \textit{center} $\zeta^{\textrm{left}}(L)$ (respectively $\zeta^{\textrm{right}}(L)$) of a Leibniz algebra $L$ is defined by the rule:
$$\zeta^{\mathrm{left}}(L)=\{x\in L|\ [x,y]=0\ \mbox{for each element }y\in L\}$$
(respectively,
$$\zeta^{\mathrm{right}}(L)=\{x\in L|\ [y,x]=0\ \mbox{for each element }y\in L\}).$$
It is not hard to prove that the left center of $L$ is an ideal, but that is not true for the right center. Moreover, $\mathrm{Leib}(L)\leqslant\zeta^{\mathrm{left}}(L)$ so that $L/\zeta^{\mathrm{left}}(L)$ is a Lie algebra. The right center is a subalgebra of $L$ and, in general, the left and right centers are different (see, for example, \cite{KOP2016}).

The \textit{center} $\zeta(L)$ of $L$ is defined by the rule:
$$\zeta(L)=\{x\in L|\ [x,y]=0=[y,x]\ \mbox{for each element }y\in L\}.$$
The center is an ideal of $L$.

Now we define the \textit{upper central series}
$$\langle0\rangle=\zeta_{0}(L)\leqslant\zeta_{1}(L)\leqslant\ldots\zeta_{\alpha}(L)\leqslant\zeta_{\alpha+1}(L)\leqslant\ldots\zeta_{\eta}(L)=\zeta_{\infty}(L)$$
of a Leibniz algebra $L$ by the following rule: $\zeta_{1}(L)=\zeta(L)$ is the center of $L$, and recursively, $\zeta_{\alpha+1}(L)/\zeta_{\alpha}(L)=\zeta(L/\zeta_{\alpha}(L))$ for all ordinals $\alpha$, and $\zeta_{\lambda}(L)=\bigcup_{\mu<\lambda}\zeta_{\mu}(L)$ for the limit ordinals $\lambda$. By definition, each term of this series is an ideal of $L$.

Define the \textit{lower central series} of $L$
$$L=\gamma_{1}(L)\geqslant\gamma_{2}(L)\geqslant\ldots\gamma_{\alpha}(L)\geqslant\gamma_{\alpha+1}\geqslant\ldots\gamma_{\tau}(L)=\gamma_{\infty}(L)$$
by the rule: $\gamma_{1}(L)=L$, $\gamma_{2}(L)=[L,L]$, and recursively $\gamma_{\alpha+1}(L)=[L,\gamma_{\alpha}(L)]$ for all ordinals $\alpha$ and $\gamma_{\lambda}(L)=\bigcap_{\mu<\lambda}\gamma_{\mu}(L)$ for the limit ordinals $\lambda$.

We say that a Leibniz algebra $L$ is \textit{nilpotent}, if there exists a positive integer $k$ such that $\gamma_{k}(L)=\langle0\rangle$. More precisely, $L$ is said to be \textit{nilpotent of nilpotency class $c$} if $\gamma_{c+1}(L)=\langle0\rangle$, but $\gamma_{c}(L)\neq\langle0\rangle$. We denote the nilpotency class of $L$ by $ncl(L)$.

Let $L$ be a Leibniz algebra over a field $F$, $M$ be non-empty subset of $L$ and $H$ be a subalgebra of $L$. Put
\begin{gather*}
\mathrm{Ann}_{H}^{\mathrm{left}}(M)=\{a\in H|\ [a,M]=\langle0\rangle\},\\
\mathrm{Ann}_{H}^{\mathrm{right}}(M)=\{a\in H|\ [M,a]=\langle0\rangle\}.
\end{gather*}
The subset $\mathrm{Ann}_{H}^{\mathrm{left}}(M)$ is called the \textit{left annihilator} of $M$ in subalgebra $H$. The subset $\mathrm{Ann}_{H}^{\mathrm{right}}(M)$ is called the \textit{right annihilator} of $M$ in subalgebra $H$. The intersection
\begin{gather*}
\mathrm{Ann}_{H}(M)=\mathrm{Ann}_{H}^{\mathrm{left}}(M)\cap\mathrm{Ann}_{H}^{\mathrm{right}}(M)=\\
\{a\in H|\ [a,M]=\langle0\rangle=[M,a]\}
\end{gather*}
is called the \textit{annihilator} of $M$ in subalgebra $H$. It is not hard to see that all of these subsets are subalgebras of $L$. Moreover, if $M$ is an ideal of $L$, then $\mathrm{Ann}_{L}(M)$ is an ideal of $L$ (see, for example, \cite{KKPS2017}).

Let $L$ be a nilpotent Leibniz algebra, whose nilpotency class is 2 and the center of $L$ has dimension 1. Of course we will suppose that $L$ is not a Lie algebra. Then there is an element $a_{1}$ such that $[a_{1},a_{1}]=a_{3}\neq0$. Since $L/\zeta(L)$ is abelian, $a_{3}\in\zeta(L)$. It follows that $[a_{1},a_{3}]=[a_{3},a_{1}]=[a_{3},a_{3}]=0$. Then $\zeta(L)=Fa_{3}$. For every element $x\in L$ we have: $[a_{1},x],[x,a_{1}]\in\zeta(L)\leqslant\langle a_{1}\rangle=Fa_{1}\oplus Fa_{3}$. It follows that a subalgebra $\langle a_{1}\rangle$ is an ideal of $L$. Since $\mathrm{dim}_{F}(\langle a_{1}\rangle)=2$, $\langle a_{1}\rangle\neq L$. Choose an element $b$ such that $b\not\in\langle a_{1}\rangle$. Then $[b,a_{1}]=\gamma a_{3}$ for some $\gamma\in F$. If $\gamma\neq0$, then put $b_{1}=\gamma^{-1}b-a_{1}$. Then $[b_{1},a_{1}]=0$. The choice of $b_{1}$ shows that $b_{1}\not\in\langle a_{1}\rangle$. If follows that a subalgebra $\mathrm{Ann}^{\mathrm{left}}_{L}(a_{1})$ has dimension 2. The paper \cite{KPSe2023} considered the case when $\mathrm{Ann}^{\mathrm{left}}_{L}(a_{1})$ is an abelian subalgebra. The next natural step is to study the situation when $\mathrm{Ann}^{\mathrm{left}}_{L}(a_{1})$ is a non-abelian subalgebra. In this case $[x,x]\neq0$ for each element $x\in\mathrm{Ann}^{\mathrm{left}}_{L}(a_{1})$, where $x\not\in\zeta(L)$. It follows that a subalgebra $\mathrm{Ann}^{\mathrm{left}}_{L}(a_{1})$ is a one-generator nilpotent algebra of dimension 2. Moreover, $\mathrm{Ann}^{\mathrm{left}}_{L}(a_{1})$ is an ideal of $L$, because
$$[L,L]=\zeta(L)\leqslant\mathrm{Ann}^{\mathrm{left}}_{L}(a_{1}).$$
Let $b$ be an element, generated $\mathrm{Ann}^{\mathrm{left}}_{L}(a_{1})$. Since $\mathrm{Ann}^{\mathrm{left}}_{L}(a_{1})$ is non-abelian, $b\not\in\zeta(L)$. We have $[b,a_{1}]=0$ and $[a_{1},b]=\gamma a_{3}$ for some $\gamma\in F$. If $\gamma=0$, then the fact that $\mathrm{Ann}^{\mathrm{left}}_{L}(a_{1})=\langle b\rangle$ shows that $\mathrm{Ann}^{\mathrm{left}}_{L}(a_{1})=\mathrm{Ann}_{L}(a_{1})$, so that $[\langle a_{1}\rangle,\langle b\rangle]=\langle0\rangle$. Thus, we obtain the following type of nilpotent Leibniz algebras:
\begin{gather*}
\mathrm{Lei}_{4}(3,F)=Fa_{1}\oplus Fa_{2}\oplus Fa_{3},\ \mbox{where }[a_{1},a_{1}]=a_{3},\\
[a_{2},a_{2}]=\lambda a_{3}, 0\neq\lambda\in F,\\
[a_{1},a_{2}]=[a_{1},a_{3}]=[a_{2},a_{1}]=[a_{2},a_{3}]=[a_{3},a_{1}]=[a_{3},a_{2}]=[a_{3},a_{3}]=0.
\end{gather*}
In other words, $\mathrm{Lei}_{4}(3,F)=L$ is the sum of two ideals $A_{1}=Fa_{1}\oplus Fa_{3}$ and $A_{2}=Fa_{2}\oplus Fa_{3}$, where $A_{1},A_{2}$ are nilpotent cyclic Leibniz algebras of dimension 2, $[A_{1},A_{2}]=[A_{2},A_{1}]=\langle0\rangle$, $\mathrm{Leib}(L)=[L,L]=\zeta^{\mathrm{left}}(L)=\zeta^{\mathrm{right}}(L)=\zeta(L)=Fa_{3}$.

This article is devoted to the description of this type of nilpotent Leibniz algebras.

\section{The description of the automorphism group of Leibniz algebras of type $\mathrm{Lei}_{4}(3,F)$.}
Let $x$ be an arbitrary element of $\mathrm{Lei}_{4}(3,F)$, $x=\xi_{1}a_{1}+\xi_{2}a_{2}+\xi_{3}a_{3}$. We have
\begin{gather*}
[x,x]=[\xi_{1}a_{1}+\xi_{2}a_{2}+\xi_{3}a_{3},\xi_{1}a_{1}+\xi_{2}a_{2}+\xi_{3}a_{3}]=\\
\xi_{1}^{2}[a_{1},a_{1}]+\xi_{2}^{2}[a_{2},a_{2}]=(\xi_{1}^{2}+\lambda\xi_{2}^{2})a_{3}.
\end{gather*}
If we suppose that $[x,x]=0$ then we come to the Leibniz algebras, whose automorphism groups have been considered in the papers \cite{KPSe2023,KPV2023}. Therefore we will suppose that $[x,x]\neq0$. If $\xi_{1}=0$ or $\xi_{2}=0$, then $[x,x]\neq0$. Suppose that $\xi_{1}\neq0$ and $\xi_{2}\neq0$. Then we can see that a polynomial $X^{2}+\lambda$ has no root in a field $F$.

We say that a field $F$ is \textit{$2$-closed}, if the equation $X^{2}=a$ has a solution in $F$ for every element $a\neq0$.

Note that every locally finite (in particular, finite) field of characteristic 2 is 2-closed. Thus Leibniz algebras of the type $\mathrm{Lei}_{4}(3,F)$ over a 2-closed field $F$ may not exist.

The Leibniz algebra $L$ is called \textit{extraspecial} if $[L,L]=\zeta(L)$ is an ideal of dimension 1. Thus, the Leibniz algebra of type $\mathrm{Lei}_{4}(3,F)$ is extraspecial.

We present some general properties of endomorphisms, automorphisms and automorphisms groups of Leibniz algebras, proofs of which can be found in \cite{KPV2023}.

\begin{lem}\label{L1}
Let $L$ be a Leibniz algebra over a field $F$, $f$ be an automorphism of $L$. Then $f(\zeta^{\mathrm{left}}(L))=\zeta^{\mathrm{left}}(L)$, $f(\zeta^{\mathrm{right}}(L))=\zeta^{\mathrm{right}}(L)$, $f(\zeta(L))=\zeta(L)$, $f([L,L])=[L,L]$.
\end{lem}

\begin{lem}\label{L2}
Let $L$ be a Leibniz algebra over a field $F$, $f$ be an automorphism of $L$. Then $f(\zeta_{\alpha}(L))=\zeta_{\alpha}(L)$, $f(\gamma_{\alpha}(L))=\gamma_{\alpha}(L)$ for all ordinals $\alpha$. In particular, $f(\zeta_{\infty}(L))=\zeta_{\infty}(L)$ and $f(\gamma_{\infty}(L))=\gamma_{\infty}(L)$.
\end{lem}

\begin{lem}\label{L3}
Let $L$ be a Leibniz algebra over a field $F$, $f$ be an endomorphism of $L$. Then $f(\gamma_{\alpha}(L))\leqslant\gamma_{\alpha}(L)$ for all ordinals $\alpha$. In particular, $f(\gamma_{\infty}(L))\leqslant\gamma_{\infty}(L)$.
\end{lem}

Let $L$ be a Leibniz algebra over a field $F$, $A$ be a subalgebra of $L$, $G=\mathrm{Aut}_{[,]}(L)$. Then we put
$$C_{G}(A)=\{\alpha\in G|\ \alpha(x)=x\ \mbox{for every }x\in A\}.$$
If $A$ is an ideal of $L$, then we put
\begin{gather*}
C_{G}(L/A)=\{\alpha\in G|\ \alpha(x+A)=x+A\ \mbox{for every }x\in L\}=\\
\{\alpha\in G|\ \alpha(x)\in x+A\ \mbox{for every }x\in L\}.
\end{gather*}

\begin{lem}\label{L4}
Let $L$ be a Leibniz algebra over a field $F$, $G=\mathrm{Aut}[,](L)$. If $A$ is a $G$-invariant subalgebra of $L$, then $C_{G}(A)$ and $C_{G}(L/A)$ are normal subgroup of $G$.
\end{lem}

Consider the automorphism groups of extraspecial Leibniz algebras.

Let $V$ be a vector space over a field $F$ and suppose that $\Phi$ is a bilinear form on $V$. We say that an automorphism $f$ of a vector space $V$ \textit{preserve} a bilinear form $\Phi$ if $\Phi(f(x),f(y))=\Phi(x,y)$ for all elements $x,y\in V$.

Denote by $\mathrm{B}(V,\Phi)$ the subset of automorphisms of a vector space $V$, preserving a bilinear form $\Phi$. It is clear that $\mathrm{B}(V,\Phi)$ is a subgroup of $\mathrm{GL}(V,F)$.

The structure of the automorphism groups of vector spaces preserving bilinear form was studied, for example, in \cite{Sz2005,Sz2013}.

\begin{lem}\label{L5}
Let $L$ be an extraspecial Leibniz algebra over a field $F$, $Z=\zeta(L)=Fc$, $V=L/Z$, $G=\mathrm{Aut}_{[,]}(L)$. Define the mapping $\Phi:V\times V\rightarrow F$ by the rule $\Phi(x+Z,y+Z)=\sigma_{xy}$, where $[x,y]=\sigma_{xy}c$. Then $\Phi$ is a bilinear form and $G/C_{G}(L/Z)$ is isomorphic to some subgroup of the automorphism group of a vector space $V$ preserving the bilinear form $\Phi$.
\end{lem}
\pf Let's fix the element $c$. Let $x+Z,y+Z$ are arbitrary cosets. Then $[x,y]=\sigma_{xy}c$, where $\sigma_{xy}\in F$. If $x_{1},y_{1}$ are elements of $L$ such that $x_{1}+Z=x+Z$, $y_{1}+Z=y+Z$, then $x_{1}=x+z_{1}$, $y_{1}=y+z_{2}$ for some elements $z_{1},z_{2}\in Z$. We have $[x_{1},y_{1}]=[x+z_{1},y+z_{2}]=[x,y]$. Thus we can see that a definition of the form $\Phi$ is correct. The fact that the operation $[,]$ is bilinear means that the form $\Phi$ is bilinear.

Let $f\in\mathrm{Aut}_{[,]}(L)$. Using Lemma~\ref{L1} we obtain that $f(Z)=Z$. Define the mapping $f^{\uparrow}:L/Z\rightarrow L/Z$ by the rule $f^{\uparrow}(x+Z)=f(x)+Z$. It is possible to prove that $f^{\uparrow}$ is a linear transformation of a vector space $V$. Moreover, this transformation is non-degenerate. We have $[x,y]=\sigma_{xy}c$, $[f(x),f(y)]=\sigma_{f(x)f(y)}c$. Since $[f(x),f(y)]=[x,y]$, then $\sigma_{f(x)f(y)}=\sigma_{xy}$. Thus
\begin{gather*}
\Phi(f^{\uparrow}(x+Z),f^{\uparrow}(y+Z))=\Phi(f(x)+Z,f(y)+Z)=\\
\sigma_{f(x)f(y)}=\sigma_{xy}=\Phi(x+Z,y+Z).
\end{gather*}
This equality shows that $f^{\uparrow}$ is an automorphism of a vector space $V$ preserving a bilinear form $\Phi$.

Consider now the mapping $\psi:G\rightarrow\mathrm{B}(V,\phi)$, defined by the rule $\psi(f)=f^{\uparrow}$, $f\in G$. It is not hard to prove that $\psi$ is a homomorphism and $\mathrm{Ker}(\psi)=C_{G}(L/Z)$. By Lemma~\ref{L4} a subgroup $C_{G}(L/Z)$ is normal in $G$. Thus we can see that a factor-group $G/C_{G}(L/Z)$ is isomorphic to a subgroup of $\mathrm{B}(V,\Phi)$. \qed

\begin{lem}\label{L6}
Let $L$ be an extraspecial Leibniz algebra over a field $F$, $Z=\zeta(L)$, $\mathrm{dim}_{F}(L)=n+1$, $G=\mathrm{Aut}_{[,]}(L)$. Then $C_{G}(L/Z)$ is isomorphic to a subgroup of $\mathrm{GL}_{n+1}(F)$, which consists of matrices of the following form:
\begin{equation*}
\left(\begin{array}{cccccc}
1 & 0 & 0 & \ldots & 0 & 0\\
0 & 1 & 0 & \ldots & 0 & 0\\
0 & 0 & 1 & \ldots & 0 & 0\\
\ldots & \ldots & \ldots & \ldots & \ldots & \ldots\\
\alpha_{1} & \alpha_{2} & \alpha_{3} & \ldots & \alpha_{n} & 1
\end{array}\right),
\end{equation*}
$\alpha_{j}\in F$, $1\leqslant j\leqslant n$. In particular, $C_{G}(L/Z)$ is isomorphic to direct product of $n$ copies of additive group of a field $F$.
\end{lem}
\pf Let $\{a_{1},a_{2},\ldots,a_{n},c\}$ be a basis of $L$. If $f\in C_{G}(L/Z)$, then $f(a_{j})=a_{j}+\alpha_{j}c$, $1\leqslant j\leqslant n$. Denote by $\Xi$ a canonical monomorphism of $C_{G}(L/Z)$ in $\mathrm{GL}_{n+1}(F)$. Then $\Xi(C_{G}(L/Z))$ is a subgroup of $\mathrm{GL}_{n+1}(F)$, which consists of matrices of the following form:
\begin{equation*}
\left(\begin{array}{cccccc}
1 & 0 & 0 & \ldots & 0 & 0\\
0 & 1 & 0 & \ldots & 0 & 0\\
0 & 0 & 1 & \ldots & 0 & 0\\
\ldots & \ldots & \ldots & \ldots & \ldots & \ldots\\
\alpha_{1} & \alpha_{2} & \alpha_{3} & \ldots & \alpha_{n} & 1
\end{array}\right),
\end{equation*}
$\alpha_{j}\in F$, $1\leqslant j\leqslant n$. It is not hard to see, that this subgroup (and hence $C_{G}(L/Z)$) is isomorphic to direct product of $n$ copies of additive group of a field $F$. \qed

\begin{thm}\label{T1}
Let $G$ be an automorphism group of Leibniz algebra $\mathrm{Lei}_{4}(3,F)$.

If $\mathrm{char}(F)=2$, then $G$ is isomorphic to a subgroup of $\mathrm{GL}_{3}(F)$, which consists of matrices of the following form:
\begin{equation*}
\left(\begin{array}{ccc}
\alpha_{1} & \lambda\alpha_{2} & 0\\
\alpha_{2} & \alpha_{1} & 0\\
\alpha_{3} & \beta_{3} & \alpha_{1}^{2}+\lambda\alpha_{2}^{2}
\end{array}\right),
\end{equation*}
$\alpha_{1},\alpha_{2},\alpha_{3},\beta_{3}\in F$. Furthermore, $G$ has a normal subgroup $C=C_{G}(L/\zeta(L))$, which is isomorphic to direct product of two copies of additive group of a field $F$ and $G/C$ is isomorphic to a subgroup of $\mathrm{GL}_{2}(F)$, which consists of matrices of the following form:
\begin{equation*}
\left(\begin{array}{cc}
\alpha_{1} & \lambda\alpha_{2}\\
\alpha_{2} & \alpha_{1}\\
\end{array}\right).
\end{equation*}

If $\mathrm{char}(F)\neq2$, then $G$ is isomorphic to a subgroup of $\mathrm{GL}_{3}(F)$, which consists of matrices of the following form:
\begin{equation*}
\left(\begin{array}{ccc}
\alpha_{1} & \delta\lambda\alpha_{2} & 0\\
\alpha_{2} & -\delta\alpha_{1} & 0\\
\alpha_{3} & \beta_{3} & \alpha_{1}^{2}+\lambda\alpha_{2}^{2}
\end{array}\right),
\end{equation*}
$\alpha_{1},\alpha_{2},\alpha_{3},\beta_{3}\in F$, $\delta\in\{-1,1\}$. Furthermore, $G$ has a normal subgroup $C=C_{G}(L/\zeta(L))$, which is isomorphic to direct product of two copies of additive group of a field $F$ and $G/C$ is isomorphic to a subgroup of $\mathrm{GL}_{2}(F)$, which consists of matrices of the following form:
\begin{equation*}
\left(\begin{array}{cc}
\alpha_{1} & \delta\lambda\alpha_{2}\\
\alpha_{2} & -\delta\alpha_{1}\\
\end{array}\right).
\end{equation*}
\end{thm}
\pf Let $L=\mathrm{Lei}_{4}(3,F)$, $f\in\mathrm{Aut}_{[,]}(L)$. By Lemma~\ref{L1} $f(Fa_{3})=Fa_{3}$. We have
\begin{gather*}
f(a_{1})=\alpha_{1}a_{1}+\alpha_{2}a_{2}+\alpha_{3}a_{3},\\
f(a_{2})=\beta_{1}a_{1}+\beta_{2}a_{2}+\beta_{3}a_{3}.
\end{gather*}
Then
\begin{gather*}
f(a_{3})=f([a_{1},a_{1}])=[f(a_{1}),f(a_{1})]=[\alpha_{1}a_{1}+\alpha_{2}a_{2}+\alpha_{3}a_{3},\alpha_{1}a_{1}+\alpha_{2}a_{2}+\alpha_{3}a_{3}]=\\
\alpha_{1}^{2}[a_{1},a_{1}]+\alpha_{2}^{2}[a_{2},a_{2}]=\alpha_{1}^{2}a_{3}+\lambda\alpha_{2}^{2}a_{3}=(\alpha_{1}^{2}+\lambda\alpha_{2}^{2})a_{3};\\
f(a_{3})=\lambda^{-1}f([a_{2},a_{2}])=\lambda^{-1}[f(a_{2}),f(a_{2})]=\\
\lambda^{-1}[\beta_{1}a_{1}+\beta_{2}a_{2}+\beta_{3}a_{3},\beta_{1}a_{1}+\beta_{2}a_{2}+\beta_{3}a_{3}]=\\
\lambda^{-1}\beta_{1}^{2}[a_{1},a_{1}]+\lambda^{-1}\beta_{2}^{2}[a_{2},a_{2}]=\lambda^{-1}\beta_{1}^{2}a_{3}+\lambda^{-1}\lambda\beta_{2}^{2}a_{3}=(\lambda^{-1}\beta_{1}^{2}+\beta_{2}^{2})a_{3};\\
0=f([a_{1},a_{2}])=[f(a_{1}),f(a_{2})]=[\alpha_{1}a_{1}+\alpha_{2}a_{2}+\alpha_{3}a_{3},\beta_{1}a_{1}+\beta_{2}a_{2}+\beta_{3}a_{3}]=\\
\alpha_{1}\beta_{1}[a_{1},a_{1}]+\alpha_{2}\beta_{2}[a_{2},a_{2}]=\alpha_{1}\beta_{1}a_{3}+\lambda\alpha_{2}\beta_{2}a_{3}=(\alpha_{1}\beta_{1}+\lambda\alpha_{2}\beta_{2})a_{3}.
\end{gather*}
Thus $\alpha_{1}^{2}+\lambda\alpha_{2}^{2}=\lambda^{-1}\beta_{1}^{2}+\beta_{2}^{2}$, $\alpha_{1}\beta_{1}+\lambda\alpha_{2}\beta_{2}=0$.

Denote by $\Xi$ a canonical monomorphism of $\mathrm{Aut}_{[,]}(L)$ in $\mathrm{GL}_{3}(F)$. Then $\Xi(f)$ is a matrix of the following form:
\begin{equation*}
\left(\begin{array}{ccc}
\alpha_{1} & \beta_{1} & 0\\
\alpha_{2} & \beta_{2} & 0\\
\alpha_{3} & \beta_{3} & \alpha_{1}^{2}+\lambda\alpha_{2}^{2}
\end{array}\right),
\end{equation*}
$\alpha_{1},\alpha_{2},\alpha_{3},\beta_{1},\beta_{2}\beta_{3}\in F$, where $\alpha_{1}^{2}+\lambda\alpha_{2}^{2}=\lambda^{-1}\beta_{1}^{2}+\beta_{2}^{2}$, $\alpha_{1}\beta_{1}+\lambda\alpha_{2}\beta_{2}=0$. In particular, if $\lambda=1$, then $\alpha_{1}^{2}+\alpha_{2}^{2}=\beta_{1}^{2}+\beta_{2}^{2}$, $\alpha_{1}\beta_{1}+\alpha_{2}\beta_{2}=0$.

Conversely, let $f$ be a linear transformation of $L$, which in a basis $\{a_{1},a_{2},a_{3}\}$ has the above matrix. Let $x,y$ be the arbitrary elements of $L$, $x=\xi_{1}a_{1}+\xi_{2}a_{2}+\xi_{3}a_{3}$, $y=\eta_{1}a_{1}+\eta_{2}a_{2}+\eta_{3}a_{3}$, where $\xi_{1},\xi_{2},\xi_{3},\eta_{1},\eta_{2},\eta_{3}\in F$. Then
\begin{gather*}
[x,y]=[\xi_{1}a_{1}+\xi_{2}a_{2}+\xi_{3}a_{3},\eta_{1}a_{1}+\eta_{2}a_{2}+\eta_{3}a_{3}]=\xi_{1}\eta_{1}[a_{1},a_{1}]+\xi_{2}\eta_{2}[a_{2},a_{2}]=\\
\xi_{1}\eta_{1}a_{3}+\lambda\xi_{2}\eta_{2}a_{3}=(\xi_{1}\eta_{1}+\lambda\xi_{2}\eta_{2})a_{3};\\
f(x)=f(\xi_{1}a_{1}+\xi_{2}a_{2}+\xi_{3}a_{3})=\xi_{1}f(a_{1})+\xi_{2}f(a_{2})+\xi_{3}f(a_{3})=\\
\xi_{1}(\alpha_{1}a_{1}+\alpha_{2}a_{2}+\alpha_{3}a_{3})+\xi_{2}(\beta_{1}a_{1}+\beta_{2}a_{2}+\beta_{3}a_{3})+\xi_{3}(\alpha_{1}^{2}+\lambda\alpha_{2}^{2})a_{3}=\\
(\xi_{1}\alpha_{1}+\xi_{2}\beta_{1})a_{1}+(\xi_{1}\alpha_{2}+\xi_{2}\beta_{2})a_{2}+(\xi_{1}\alpha_{3}+\xi_{2}\beta_{3}+\xi_{3}\alpha_{1}^{2}+\xi_{3}\lambda\alpha_{2}^{2})a_{3};\\
f(y)=(\eta_{1}\alpha_{1}+\eta_{2}\beta_{1})a_{1}+(\eta_{1}\alpha_{2}+\eta_{2}\beta_{2})a_{2}+(\eta_{1}\alpha_{3}+\eta_{2}\beta_{3}+\eta_{3}\alpha_{1}^{2}+\eta_{3}\lambda\alpha_{2}^{2})a_{3};\\
f([x,y])=f((\xi_{1}\eta_{1}+\lambda\xi_{2}\eta_{2})a_{3})=(\xi_{1}\eta_{1}+\lambda\xi_{2}\eta_{2})f(a_{3})=\\
(\xi_{1}\eta_{1}+\lambda\xi_{2}\eta_{2})(\alpha_{1}^{2}+\lambda\alpha_{2}^{2})a_{3}=(\xi_{1}\eta_{1}\alpha_{1}^{2}+\lambda\xi_{2}\eta_{2}\alpha_{1}^{2}+\xi_{1}\eta_{1}\lambda\alpha_{2}^{2}+\lambda^{2}\xi_{2}\eta_{2}\alpha_{2}^{2})a_{3};\\
[f(x),f(y)]=\\
[(\xi_{1}\alpha_{1}+\xi_{2}\beta_{1})a_{1}+(\xi_{1}\alpha_{2}+\xi_{2}\beta_{2})a_{2}+(\xi_{1}\alpha_{3}+\xi_{2}\beta_{3}+\xi_{3}\alpha_{1}^{2}+\xi_{3}\lambda\alpha_{2}^{2})a_{3},\\
(\eta_{1}\alpha_{1}+\eta_{2}\beta_{1})a_{1}+(\eta_{1}\alpha_{2}+\eta_{2}\beta_{2})a_{2}+(\eta_{1}\alpha_{3}+\eta_{2}\beta_{3}+\eta_{3}\alpha_{1}^{2}+\eta_{3}\lambda\alpha_{2}^{2})a_{3}]=\\
(\xi_{1}\alpha_{1}+\xi_{2}\beta_{1})(\eta_{1}\alpha_{1}+\eta_{2}\beta_{1})[a_{1},a_{1}]+(\xi_{1}\alpha_{2}+\xi_{2}\beta_{2})(\eta_{1}\alpha_{2}+\eta_{2}\beta_{2})[a_{2},a_{2}]=\\
(\xi_{1}\alpha_{1}+\xi_{2}\beta_{1})(\eta_{1}\alpha_{1}+\eta_{2}\beta_{1})a_{3}+\lambda(\xi_{1}\alpha_{2}+\xi_{2}\beta_{2})(\eta_{1}\alpha_{2}+\eta_{2}\beta_{2})a_{3}=\\
(\xi_{1}\alpha_{1}\eta_{1}\alpha_{1}+\xi_{1}\alpha_{1}\eta_{2}\beta_{1}+\xi_{2}\beta_{1}\eta_{1}\alpha_{1}+\xi_{2}\beta_{1}\eta_{2}\beta_{1}+\\
\lambda\xi_{1}\alpha_{2}\eta_{1}\alpha_{2}+\lambda\xi_{1}\alpha_{2}\eta_{2}\beta_{2}+\lambda\xi_{2}\beta_{2}\eta_{1}\alpha_{2}+\lambda\xi_{2}\beta_{2}\eta_{2}\beta_{2})a_{3}=\\
(\xi_{1}\eta_{1}\alpha_{1}^{2}+\xi_{1}\eta_{2}\alpha_{1}\beta_{1}+\xi_{2}\eta_{1}\alpha_{1}\beta_{1}+\xi_{2}\eta_{2}\beta_{1}^{2}+\\
\lambda\xi_{1}\eta_{1}\alpha_{2}^{2}+\lambda\xi_{1}\eta_{2}\alpha_{2}\beta_{2}+\lambda\xi_{2}\eta_{1}\alpha_{2}\beta_{2}+\lambda\xi_{2}\eta_{2}\beta_{2}^{2})a_{3}=\\
(\xi_{1}\eta_{1}(\alpha_{1}^{2}+\lambda\alpha_{2}^{2})+\xi_{2}\eta_{2}(\beta_{1}^{2}+\lambda\beta_{2}^{2})+\xi_{1}\eta_{2}(\alpha_{1}\beta_{1}+\lambda\alpha_{2}\beta_{2})+\xi_{2}\eta_{1}(\alpha_{1}\beta_{1}+\lambda\alpha_{2}\beta_{2}))a_{3}.
\end{gather*}
Using an equality $f([x,y])=[f(x),f(y)]$ we obtain
\begin{gather*}
\xi_{1}\eta_{1}\alpha_{1}^{2}+\lambda\xi_{2}\eta_{2}\alpha_{1}^{2}+\xi_{1}\eta_{1}\lambda\alpha_{2}^{2}+\lambda^{2}\xi_{2}\eta_{2}\alpha_{2}^{2}=\\
\xi_{1}\eta_{1}(\alpha_{1}^{2}+\lambda\alpha_{2}^{2})+\xi_{2}\eta_{2}(\beta_{1}^{2}+\lambda\beta_{2}^{2})+\xi_{1}\eta_{2}(\alpha_{1}\beta_{1}+\lambda\alpha_{2}\beta_{2})+\xi_{2}\eta_{1}(\alpha_{1}\beta_{1}+\lambda\alpha_{2}\beta_{2}).
\end{gather*}
It follows that
$$\xi_{2}\eta_{2}(\lambda\alpha_{1}^{2}+\lambda^{2}\alpha_{2}^{2}-\beta_{1}^{2}-\lambda\beta_{2}^{2})-\xi_{1}\beta_{2}(\alpha_{1}\beta_{1}+\lambda\alpha_{2}\beta_{2})-\xi_{2}\eta_{1}(\alpha_{1}\beta_{1}+\lambda\alpha_{2}\beta_{2})=0.$$
Taking into account the equalities $\alpha_{1}^{2}+\lambda\alpha_{2}^{2}=\lambda^{-1}\beta_{1}^{2}+\beta_{2}^{2}$ and $\alpha_{1}\beta_{1}+\lambda\alpha_{2}\beta_{2}=0$, we obtain that $f([x,y])=[f(x),f(y)]$.

Suppose first that $\mathrm{char}(F)=2$. If $\alpha_{2}=0$, then $\alpha_{1}\beta_{1}=0$. In this case either $\alpha_{1}=0$ or $\beta_{1}=0$. The case $\alpha_{1}=0$ is impossible (otherwise a matrix of linear transformation $f$ is degenerate). Hence $\beta_{1}=0$. It follows that $\alpha_{1}^{2}=\beta_{2}^{2}$. Since $\mathrm{char}(F)=2$, then $\beta_{2}=\alpha_{1}$. In this case $\Xi(f)$ is a matrix of the following form:
\begin{equation*}
\left(\begin{array}{ccc}
\alpha_{1} & 0 & 0\\
0 & \alpha_{1} & 0\\
\alpha_{3} & \beta_{3} & \alpha_{1}^{2}
\end{array}\right),
\end{equation*}
$\alpha_{1},\alpha_{3},\beta_{3}\in F$.

If $\alpha_{1}=0$, then $\lambda\alpha_{2}\beta_{2}=0$. Since $\lambda\neq0$, then either $\alpha_{2}=0$ or $\beta_{2}=0$. If we suppose that $\alpha_{2}=0$, then a matrix of linear transformation $f$ is degenerate, and we obtain a contradiction. Hence $\beta_{2}=0$. It follows that $\lambda\alpha_{2}^{2}=\lambda^{-1}\beta_{1}^{2}$ that is $\beta_{1}^{2}=\lambda^{2}\alpha_{2}^{2}$. The fact that $\mathrm{char}(F)=2$ implies that $\beta_{1}=\lambda\alpha_{2}$. In this case $\Xi(f)$ is a matrix of the following form:
\begin{equation*}
\left(\begin{array}{ccc}
0 & \lambda\alpha_{2} & 0\\
\alpha_{2} & 0 & 0\\
\alpha_{3} & \beta_{3} & \lambda\alpha_{2}^{2}
\end{array}\right),
\end{equation*}
$\alpha_{2},\alpha_{3},\beta_{3}\in F$.

Suppose now that all coefficients $\alpha_{1},\alpha_{2},\beta_{1},\beta_{2}$ are non-zero. The equality $\alpha_{1}\beta_{1}+\lambda\alpha_{2}\beta_{2}=0$ implies that $\alpha_{1}\beta_{1}=\lambda\alpha_{2}\beta_{2}$, that is $\alpha_{1}\alpha_{2}^{-1}=\lambda\beta_{2}\beta_{1}^{-1}=\kappa$. Then $\alpha_{1}=\alpha_{2}\kappa$, $\beta_{2}=\lambda^{-1}\beta_{1}\kappa$. We have
$$\alpha_{1}^{2}+\lambda\alpha_{2}^{2}=\alpha_{2}^{2}\kappa^{2}+\lambda\alpha_{2}^{2}=\lambda^{-1}\beta_{1}^{2}+\beta_{2}^{2}=\lambda^{-1}\beta_{1}^{2}+\lambda^{-2}\beta_{1}^{2}\kappa^{2}=\lambda^{-2}\beta_{1}^{2}(\lambda+\kappa^{2}).$$
Hence $\alpha_{2}^{2}(\kappa^{2}+\lambda)=\lambda^{-2}\beta_{1}^{2}(\lambda+\kappa^{2})$. We have $\mathrm{det}(\Xi(f))=(\alpha_{1}^{2}+\lambda\alpha_{2}^{2})(\alpha_{1}\beta_{2}-\alpha_{2}\beta_{1})$. If we suppose that $\lambda+\kappa^{2}=0$, then $\mathrm{det}(\Xi(f))=0$, which is impossible. Thus $\lambda+\kappa^{2}\neq0$. Then we obtain that $\alpha_{2}^{2}=\lambda^{-2}\beta_{1}^{2}$, that is $\beta_{1}^{2}=\lambda^{2}\alpha_{2}^{2}$. Taking into account the equalities $\alpha_{1}^{2}+\lambda\alpha_{2}^{2}=\lambda^{-1}\beta_{1}^{2}+\beta_{2}^{2}$, we obtain that $\alpha_{1}^{2}=\beta_{2}^{2}$. The fact that $\mathrm{char}(F)=2$, implies that $\beta_{1}=\lambda\alpha_{2}$ and $\beta_{2}=\alpha_{1}$. Thus $\Xi(f)$ is a matrix of the following form:
\begin{equation*}
\left(\begin{array}{ccc}
\alpha_{1} & \lambda\alpha_{2} & 0\\
\alpha_{2} & \alpha_{1} & 0\\
\alpha_{3} & \beta_{3} & \alpha_{1}^{2}+\lambda\alpha_{2}^{2}
\end{array}\right),
\end{equation*}
$\alpha_{1},\alpha_{2},\alpha_{3},\beta_{3}\in F$. We can see that for $\alpha_{1}=0$ or $\alpha_{1}=2$ we have matrices of the forms obtained above. Thus we can conclude that $\Xi(G)$ consists of non-degenerate matrices of the following form:
\begin{equation*}
\left(\begin{array}{ccc}
\alpha_{1} & \lambda\alpha_{2} & 0\\
\alpha_{2} & \alpha_{1} & 0\\
\alpha_{3} & \beta_{3} & \alpha_{1}^{2}+\lambda\alpha_{2}^{2}
\end{array}\right),
\end{equation*}
$\alpha_{1},\alpha_{2},\alpha_{3},\beta_{3}\in F$.

By Lemma~\ref{L6} a normal subgroup $\Xi(C_{G}(L/\zeta(L)))$ consists of matrices of the following form:
\begin{equation*}
\left(\begin{array}{ccc}
1 & 0 & 0\\
0 & 1 & 0\\
\alpha_{3} & \beta_{3} & 1
\end{array}\right),
\end{equation*}
$\alpha_{3},\beta_{3}\in F$. Hence $C_{G}(L/\zeta(L))$ is isomorphic to direct product of two copies of additive group of a field $F$.

Consider now the mapping $\upsilon:\Xi(G)\rightarrow\mathrm{GL}_{2}(F)$, defined by the rule
\begin{equation*}
\left(\begin{array}{ccc}
\alpha_{1} & \lambda\alpha_{2} & 0\\
\alpha_{2} & \alpha_{1} & 0\\
\alpha_{3} & \beta_{3} & \alpha_{1}^{2}+\lambda\alpha_{2}^{2}
\end{array}\right)\rightarrow
\left(\begin{array}{cc}
\alpha_{1} & \lambda\alpha_{2}\\
\alpha_{2} & \alpha_{1}
\end{array}\right),
\end{equation*}
$\alpha_{1},\alpha_{2},\alpha_{3},\beta_{3}\in F$.
We have
\begin{gather*}
\left(\begin{array}{ccc}
\alpha_{1} & \lambda\alpha_{2} & 0\\
\alpha_{2} & \alpha_{1} & 0\\
\alpha_{3} & \beta_{3} & \alpha_{1}^{2}+\lambda\alpha_{2}^{2}
\end{array}\right)
\left(\begin{array}{ccc}
\gamma_{1} & \lambda\gamma_{2} & 0\\
\gamma_{2} & \gamma_{1} & 0\\
\gamma_{3} & \sigma_{3} & \gamma_{1}^{2}+\lambda\gamma_{2}^{2}
\end{array}\right)=
\end{gather*}

\begin{footnotesize}
\[\left(\begin{array}{ccc}
\alpha_{1}\gamma_{1}+\lambda\alpha_{2}\gamma_{2} & \lambda\alpha_{1}\gamma_{2}+\lambda\alpha_{2}\gamma_{1} & 0\\
\alpha_{2}\gamma_{1}+\alpha_{1}\gamma_{2} & \lambda\alpha_{2}\gamma_{2}+\alpha_{1}\gamma_{1} & 0\\
\alpha_{3}\gamma_{1}+\beta_{3}\gamma_{2}+(\alpha_{1}^{2}+\lambda\alpha_{2}^{2})\gamma_{3} & \lambda\alpha_{3}\gamma_{2}+\beta_{3}\gamma_{1}+(\alpha_{1}^{2}+\lambda\alpha_{2}^{2})\sigma_{3} & (\alpha_{1}^{2}+\lambda\alpha_{2}^{2})(\gamma_{1}^{2}+\lambda\gamma_{2}^{2})
\end{array}\right)\]
\end{footnotesize}
and
\begin{gather*}
\left(\begin{array}{cc}
\alpha_{1} & \lambda\alpha_{2}\\
\alpha_{2} & \alpha_{1}
\end{array}\right)
\left(\begin{array}{cc}
\gamma_{1} & \lambda\gamma_{2}\\
\gamma_{2} & \gamma_{1}
\end{array}\right)=
\left(\begin{array}{cc}
\alpha_{1}\gamma_{1}+\lambda\alpha_{2}\gamma_{2} & \lambda\alpha_{1}\gamma_{2}+\lambda\alpha_{2}\gamma_{1}\\
\alpha_{2}\gamma_{1}+\alpha_{1}\gamma_{2} & \lambda\alpha_{2}\gamma_{2}+\alpha_{1}\gamma_{1}
\end{array}\right)
\end{gather*}
Thus we can see that a mapping $\upsilon$ is a homomorphism, $\mathrm{Ker}(\upsilon)$ consists of matrices of the following form:
\begin{equation*}
\left(\begin{array}{ccc}
1 & 0 & 0\\
0 & 1 & 0\\
\alpha_{3} & \beta_{3} & 1
\end{array}\right),
\end{equation*}
$\alpha_{3},\beta_{3}\in F$. Then $\mathrm{Ker}(\upsilon)=C_{G}(L/\zeta(L))$, $\mathrm{Im}(\upsilon)$ consists of matrices of the following form:
\begin{equation*}
\left(\begin{array}{cc}
\alpha_{1} & \lambda\alpha_{2}\\
\alpha_{2} & \alpha_{1}
\end{array}\right),
\end{equation*}
$\alpha_{1},\alpha_{2}\in F$.

Suppose now that $\mathrm{char}(F)\neq2$. If $\alpha_{2}=0$, then $\alpha_{1}\beta_{1}=0$. In this case either $\alpha_{1}=0$ or $\beta_{1}=0$. The case $\alpha_{1}=0$ is impossible (otherwise a matrix of linear transformation $f$ is degenerate). Hence $\beta_{1}=0$. It follows that $\alpha_{1}^{2}=\beta_{2}^{2}$. Thus $\beta_{2}=\alpha_{1}$ or $\beta_{2}=-\alpha_{1}$. In general $\beta_{2}=\delta\alpha_{1}$, where $\delta\in\{1,-1\}$, so that $\Xi(f)$ is a matrix of the following form:
\begin{equation*}
\left(\begin{array}{ccc}
\alpha_{1} & 0 & 0\\
0 & \delta\alpha_{1} & 0\\
\alpha_{3} & \beta_{3} & \alpha_{1}^{2}
\end{array}\right),
\end{equation*}
$\alpha_{1},\alpha_{3},\beta_{3}\in F$. If $\alpha_{1}=0$, then $\lambda\alpha_{2}\beta_{2}=0$. Since $\lambda\neq0$, then either $\alpha_{2}=0$ or $\beta_{2}=0$. The case $\alpha_{2}=0$ is impossible (otherwise a matrix of linear transformation $f$ is degenerate). Hence $\beta_{2}=0$. It follows that $\lambda\alpha_{2}^{2}=\lambda^{-1}\beta_{1}^{2}$, so that $\beta_{1}^{2}=\lambda^{2}\alpha_{2}^{2}$. Thus $\beta_{1}=\lambda\alpha_{2}$ or $\beta_{1}=-\lambda\alpha_{2}$. In general $\beta_{1}=\delta\lambda\alpha_{2}$, where $\delta\in\{1,-1\}$. Hence $\Xi(f)$ is a matrix of the following form:
\begin{equation*}
\left(\begin{array}{ccc}
0 & \delta\lambda\alpha_{2} & 0\\
\alpha_{2} & 0 & 0\\
\alpha_{3} & \beta_{3} & \lambda\alpha_{2}^{2}
\end{array}\right),
\end{equation*}
$\alpha_{2},\alpha_{3},\beta_{3}\in F$.

Suppose now that all coefficients $\alpha_{1},\alpha_{2},\beta_{1},\beta_{2}$ are non-zero. The equality $\alpha_{1}\beta_{1}+\lambda\alpha_{2}\beta_{2}=0$ implies that $\alpha_{1}\beta_{1}=-\lambda\alpha_{2}\beta_{2}$, that is $\alpha_{1}\alpha_{2}^{-1}=-\lambda\beta_{2}\beta_{1}^{-1}=\kappa$. Then $\alpha_{1}=\alpha_{2}\kappa$, $\beta_{2}=-\lambda^{-1}\beta_{1}\kappa$. We have
$$\alpha_{1}^{2}+\lambda\alpha_{2}^{2}=\alpha_{2}^{2}\kappa^{2}+\lambda\alpha_{2}^{2}=\lambda^{-1}\beta_{1}^{2}+\beta_{2}^{2}=\lambda^{-1}\beta_{1}^{2}+\lambda^{-2}\beta_{1}^{2}\kappa^{2}=\lambda^{-2}\beta_{1}^{2}(\lambda+\kappa^{2}).$$
Thus $\alpha_{2}^{2}(\kappa^{2}+\lambda)=\lambda^{-2}\beta_{1}^{2}(\lambda+\kappa^{2})$. We have
$$\mathrm{det}(\Xi(f))=(\alpha_{1}^{2}+\lambda\alpha_{2}^{2})(\alpha_{1}\beta_{2}-\alpha_{2}\beta_{1}).$$
If we suppose that $\lambda+\kappa^{2}=0$, then $\mathrm{det}(\Xi(f))=0$, which is impossible. Hence $\lambda+\kappa^{2}\neq0$. Then we obtain that $\alpha_{2}^{2}=\lambda^{-2}\beta_{1}^{2}$, that is $\beta_{1}^{2}=\lambda^{2}\alpha_{2}^{2}$. Taking into account the equalities $\alpha_{1}^{2}+\lambda\alpha_{2}^{2}=\lambda^{-1}\beta_{1}^{2}+\beta_{2}^{2}$, we obtain that $\alpha_{1}^{2}=\beta_{2}^{2}$. Since $\mathrm{char}(F)\neq2$, then either $\alpha_{1}=\beta_{2}$ or $\alpha_{1}=-\beta_{2}$. Suppose that $\alpha_{1}=\beta_{2}$. Then $\beta_{2}\beta_{1}+\lambda\alpha_{2}\beta_{2}=0$. It follows that $\beta_{1}+\lambda\alpha_{2}=0$, that is $\beta_{1}=-\lambda\alpha_{2}$. In this case $\Xi(f)$ is a matrix of the following form:
\begin{equation*}
\left(\begin{array}{ccc}
\alpha_{1} & -\lambda\alpha_{2} & 0\\
\alpha_{2} & \alpha_{1} & 0\\
\alpha_{3} & \beta_{3} & \alpha_{1}^{2}+\lambda\alpha_{2}^{2}
\end{array}\right),
\end{equation*}
$\alpha_{1},\alpha_{2},\alpha_{3},\beta_{3}\in F$. Suppose that $\alpha_{1}=-\beta_{2}$. Then $-\beta_{2}\beta_{1}+\lambda\alpha_{2}\beta_{2}=0$. It follows that $-\beta_{1}+\lambda\alpha_{2}=0$, that is $\beta_{1}=\lambda\alpha_{2}$. In this case $\Xi(f)$ is a matrix of the following form:
\begin{equation*}
\left(\begin{array}{ccc}
\alpha_{1} & \lambda\alpha_{2} & 0\\
\alpha_{2} & -\alpha_{1} & 0\\
\alpha_{3} & \beta_{3} & \alpha_{1}^{2}+\lambda\alpha_{2}^{2}
\end{array}\right),
\end{equation*}
$\alpha_{1},\alpha_{2},\alpha_{3},\beta_{3}\in F$.

The following equalities
\begin{equation*}
\left(\begin{array}{ccc}
\alpha_{1} & \lambda\alpha_{2} & 0\\
\alpha_{2} & -\alpha_{1} & 0\\
\alpha_{3} & \beta_{3} & \alpha_{1}^{2}+\lambda\alpha_{2}^{2}
\end{array}\right)
\left(\begin{array}{ccc}
\sigma_{1} & \lambda\sigma_{2} & 0\\
\sigma_{2} & -\sigma_{1} & 0\\
\sigma_{3} & \tau_{3} & \sigma_{1}^{2}+\lambda\sigma_{2}^{2}
\end{array}\right)=
\end{equation*}

\begin{footnotesize}
\[\left(\begin{array}{ccc}
\alpha_{1}\sigma_{1}+\lambda\alpha_{2}\sigma_{2} & \lambda\alpha_{1}\sigma_{2}-\lambda\alpha_{2}\sigma_{1} & 0\\
\alpha_{2}\sigma_{1}-\alpha_{1}\sigma_{2} & \lambda\alpha_{2}\sigma_{2}+\alpha_{1}\sigma_{1} & 0\\
\alpha_{3}\sigma_{1}+\beta_{3}\sigma_{2}+(\alpha_{1}^{2}+\lambda\alpha_{2}^{2})\sigma_{3} & \lambda\alpha_{3}\sigma_{2}-\beta_{3}\sigma_{1}+(\alpha_{1}^{2}+\lambda\alpha_{2}^{2})\tau_{3} & (\alpha_{1}^{2}+\lambda\alpha_{2}^{2})(\sigma_{1}^{2}+\lambda\sigma_{2}^{2})
\end{array}\right),\]
\end{footnotesize}

\begin{equation*}
\left(\begin{array}{ccc}
\alpha_{1} & \lambda\alpha_{2} & 0\\
\alpha_{2} & -\alpha_{1} & 0\\
\alpha_{3} & \beta_{3} & \alpha_{1}^{2}+\lambda\alpha_{2}^{2}
\end{array}\right)
\left(\begin{array}{ccc}
\sigma_{1} & -\lambda\sigma_{2} & 0\\
\sigma_{2} & \sigma_{1} & 0\\
\sigma_{3} & \tau_{3} & \sigma_{1}^{2}+\lambda\sigma_{2}^{2}
\end{array}\right)=
\end{equation*}

\begin{footnotesize}
\[\left(\begin{array}{ccc}
\alpha_{1}\sigma_{1}+\lambda\alpha_{2}\sigma_{2} & -\lambda\alpha_{1}\sigma_{2}+\lambda\alpha_{2}\sigma_{1} & 0\\
\alpha_{2}\sigma_{1}-\alpha_{1}\sigma_{2} & -\lambda\alpha_{2}\sigma_{2}-\alpha_{1}\sigma_{1} & 0\\
\alpha_{3}\sigma_{1}+\beta_{3}\sigma_{2}+(\alpha_{1}^{2}+\lambda\alpha_{2}^{2})\sigma_{3} & -\lambda\alpha_{3}\sigma_{2}+\beta_{3}\sigma_{1}+(\alpha_{1}^{2}+\lambda\alpha_{2}^{2})\tau_{3} & (\alpha_{1}^{2}+\lambda\alpha_{2}^{2})(\sigma_{1}^{2}+\lambda\sigma_{2}^{2})
\end{array}\right),\]
\end{footnotesize}

\begin{equation*}
\left(\begin{array}{ccc}
\alpha_{1} & -\lambda\alpha_{2} & 0\\
\alpha_{2} & \alpha_{1} & 0\\
\alpha_{3} & \beta_{3} & \alpha_{1}^{2}+\lambda\alpha_{2}^{2}
\end{array}\right)
\left(\begin{array}{ccc}
\sigma_{1} & \lambda\sigma_{2} & 0\\
\sigma_{2} & -\sigma_{1} & 0\\
\sigma_{3} & \tau_{3} & \sigma_{1}^{2}+\lambda\sigma_{2}^{2}
\end{array}\right)=
\end{equation*}

\begin{footnotesize}
\[\left(\begin{array}{ccc}
\alpha_{1}\sigma_{1}-\lambda\alpha_{2}\sigma_{2} & \lambda\alpha_{1}\sigma_{2}+\lambda\alpha_{2}\sigma_{1} & 0\\
\alpha_{2}\sigma_{1}+\alpha_{1}\sigma_{2} & \lambda\alpha_{2}\sigma_{2}-\alpha_{1}\sigma_{1} & 0\\
\alpha_{3}\sigma_{1}+\beta_{3}\sigma_{2}+(\alpha_{1}^{2}+\lambda\alpha_{2}^{2})\sigma_{3} & \lambda\alpha_{3}\sigma_{2}-\beta_{3}\sigma_{1}+(\alpha_{1}^{2}+\lambda\alpha_{2}^{2})\tau_{3} & (\alpha_{1}^{2}+\lambda\alpha_{2}^{2})(\sigma_{1}^{2}+\lambda\sigma_{2}^{2})
\end{array}\right),\]
\end{footnotesize}

\begin{equation*}
\left(\begin{array}{ccc}
\alpha_{1} & -\lambda\alpha_{2} & 0\\
\alpha_{2} & \alpha_{1} & 0\\
\alpha_{3} & \beta_{3} & \alpha_{1}^{2}+\lambda\alpha_{2}^{2}
\end{array}\right)
\left(\begin{array}{ccc}
\sigma_{1} & -\lambda\sigma_{2} & 0\\
\sigma_{2} & \sigma_{1} & 0\\
\sigma_{3} & \tau_{3} & \sigma_{1}^{2}+\lambda\sigma_{2}^{2}
\end{array}\right)=
\end{equation*}

\begin{footnotesize}
\[\left(\begin{array}{ccc}
\alpha_{1}\sigma_{1}-\lambda\alpha_{2}\sigma_{2} & -\lambda\alpha_{1}\sigma_{2}-\lambda\alpha_{2}\sigma_{1} & 0\\
\alpha_{2}\sigma_{1}+\alpha_{1}\sigma_{2} & -\lambda\alpha_{2}\sigma_{2}+\alpha_{1}\sigma_{1} & 0\\
\alpha_{3}\sigma_{1}+\beta_{3}\sigma_{2}+(\alpha_{1}^{2}+\lambda\alpha_{2}^{2})\sigma_{3} & -\lambda\alpha_{3}\sigma_{2}+\beta_{3}\sigma_{1}+(\alpha_{1}^{2}+\lambda\alpha_{2}^{2})\tau_{3} & (\alpha_{1}^{2}+\lambda\alpha_{2}^{2})(\sigma_{1}^{2}+\lambda\sigma_{2}^{2})
\end{array}\right),\]
\end{footnotesize}
shows that $\Xi(G)$ consists of matrices of the following form:
\begin{equation*}
\left(\begin{array}{ccc}
\alpha_{1} & \delta\lambda\alpha_{2} & 0\\
\alpha_{2} & -\delta\alpha_{1} & 0\\
\alpha_{3} & \beta_{3} & \alpha_{1}^{2}+\lambda\alpha_{2}^{2}
\end{array}\right),
\end{equation*}
$\alpha_{1},\alpha_{2},\alpha_{3},\beta_{3}\in F$, $\delta\in\{-1,1\}$.

Consider now the mapping $\upsilon:\Xi(G)\rightarrow\mathrm{GL}_{2}(F)$, defined by the rule
\begin{equation*}
\left(\begin{array}{ccc}
\alpha_{1} & \delta\lambda\alpha_{2} & 0\\
\alpha_{2} & -\delta\alpha_{1} & 0\\
\alpha_{3} & \beta_{3} & \alpha_{1}^{2}+\lambda\alpha_{2}^{2}
\end{array}\right)\rightarrow
\left(\begin{array}{cc}
\alpha_{1} & \delta\lambda\alpha_{2}\\
\alpha_{2} & -\delta\alpha_{1}\\
\end{array}\right),
\end{equation*}
$\alpha_{1},\alpha_{2},\alpha_{3},\beta_{3}\in F$, $\delta\in\{-1,1\}$.
As above we can check that this mapping is a homomorphism and $\mathrm{Ker}(\Xi)$ consists of matrices of the following form:
\begin{equation*}
\left(\begin{array}{ccc}
1 & 0 & 0\\
0 & 1 & 0\\
\alpha_{3} & \beta_{3} & 1
\end{array}\right),
\end{equation*}
$\alpha_{3},\beta_{3}\in F$.
Hence $\mathrm{Ker}(\Xi)=C_{G}(L/\zeta(L))$, $\mathrm{Im}(\Xi)$ consists of matrices of the following form:
\begin{equation*}
\left(\begin{array}{cc}
\alpha_{1} & \delta\lambda\alpha_{2}\\
\alpha_{2} & -\delta\alpha_{1}\\
\end{array}\right),
\end{equation*}
$\alpha_{1},\alpha_{2}\in F$, $\delta\in\{-1,1\}$.
The theorem is proved. \qed

\end{document}